\documentclass{article}
\usepackage{amsmath,amssymb,amsfonts,amsthm}
\usepackage{mathtools}
\usepackage{mathrsfs}
\usepackage{tasks}
\usepackage[inline]{enumitem}
\usepackage[hidelinks]{hyperref} 
\usepackage{tikz}
\usepackage{xparse}
\usepackage{xspace}
\usepackage{macros} 

\title{The finite frame property of some extensions of the pure logic of necessitation}
\author{Taishi Kurahashi\footnote{Email: kurahashi@people.kobe-u.ac.jp}
\footnote{Graduate School of System Informatics,
Kobe University,
1-1 Rokkodai, Nada, Kobe 657-8501, Japan.}
and Yuta Sato\footnote{Email: 231x032x@gsuite.kobe-u.ac.jp}
\footnote{Graduate School of System Informatics,
Kobe University,
1-1 Rokkodai, Nada, Kobe 657-8501, Japan.}}
\date{}

\begin{document}

\maketitle

\begin{abstract}
We study the finite frame property of some extensions of Fitting, Marek, and Truszczy\'nski's pure logic of necessitation $\mathbf{N}$. 
For any natural numbers $m, n$, we introduce the logic $\NRAmn$ by adding the single axiom scheme $\Box^n \varphi \to \Box^m \varphi$ and the rule $\dfrac{\neg \Box \varphi}{\neg \Box \Box \varphi}$ ($\text{\textsc{Ros}}^\Box$) into $\N$.
We prove the finite frame property of $\mathbf{N}^+\mathbf{A}_{m, n}$ with respect to Fitting, Marek, and Truszczy\'nski's relational semantics.
We also prove that for $n \ge 2$, the logic obtained by removing the rule $\text{\textsc{Ros}}^\Box$ from $\mathbf{N}^+\mathbf{A}_{0, n}$ is incomplete with respect to that semantics.
\end{abstract}

\section{Introduction}

\FMT \cite{fmt} introduced the pure logic of necessitation $\N$ with the motivation of analyzing non-monotonic reasoning (See also \cite{MT1993}). 
The logic $\N$ is a non-normal modal logic whose axioms are propositional tautologies and whose inference rules are modus ponens (\textsc{MP}) and necessitation (\textsc{Nec}). 
Namely, the logic $\N$ is obtained from classical propositional logic by simply adding the rule \textsc{Nec}. 
\FMT also introduced a natural Kripke-like semantics based on the notion of $\N$-frames and proved the finite frame property of $\N$ with respect to their semantics\footnote{A different semantics for $\N$ from that of Fitting et al.~was given and studied by Omori and Skurt \cite{OS}, where $\N$ is called $\mathbf{M^+}$.}.

Recently, the first author discovered a new aspect of the logic $\N$ in the context of provability logic \cite{Kurahashi2022}. 
Let $T$ be a recursively axiomatized consistent extension of Peano arithmetic. 
We say that a $\Sigma_1$ formula $\mathrm{Pr}_T(x)$ of arithmetic is a provability predicate of $T$ if it defines the set of all theorems of $T$ in the standard model of arithmetic. 
Then, by the $\Sigma_1$-completeness, for any sentence $\varphi$, if $\varphi$ is provable in $T$, then $\mathrm{Pr}_T(\ulcorner \varphi \urcorner)$ is also provable in $T$. 
When we interpret the modal operator $\Box$ as $\mathrm{Pr}_T(x)$, this property corresponds to the rule \textsc{Nec}. 
In \cite{Kurahashi2022}, it is proved that $T$-provable principle of $\Box$ common to all provability predicates of $T$ is only \textsc{Nec}, that is, the logic $\N$ exactly coincides with the provability logic of all provability predicates.

In the context of the second incompleteness theorem, it is important to consider provability predicates satisfying some additional conditions. 
The paper \cite{Kurahashi2022} also studied provability predicates of $T$ satisfying the condition \textsf{D3}: $T$ proves $\mathrm{Pr}_T(\ulcorner \varphi \urcorner) \to \mathrm{Pr}_T(\ulcorner \mathrm{Pr}_T(\ulcorner \varphi \urcorner) \urcorner)$ for any sentence $\varphi$. 
The logic $\mathbf{N4}$ is obtained from $\N$ by adding the axiom scheme $\Box \varphi \to \Box \Box \varphi$. 
Then, it was proved that $\mathbf{N4}$ has the finite frame property with respect to transitive $\N$-frames, and by using this property, it was also proved that $\mathbf{N4}$ is exactly the provability logic of all provability predicates satisfying \textsf{D3}. 
See \cite{Boolos1993} for details on provability predicates, incompleteness theorems, and provability logic. 

The finite frame property of $\mathbf{N4}$ with respect to transitive $\N$-frames is an analogue of the well-known fact that the normal modal logic $\mathbf{K4}$ has the finite frame property with respect to transitive Kripke frames. 
Of course, many normal modal logics other than $\mathbf{K4}$ also have finite frame property, but on the other hand, it is unknown whether the normal modal logic obtained by adding the axiom scheme $\Amn$: $\Box^n \varphi \to \Box^m \varphi$ to the logic $\mathbf{K}$ has finite frame property in general (cf.~\cite[Problem 12.1]{CZ97}, \cite[Problem 6]{WZ07}, \cite{Zak97} and \cite[Problem 6.12]{Zak97_2}). 
In this context, the following natural question arises: does the logic obtained by adding $\Amn$ to $\N$ instead of $\mathbf{K}$ have finite frame property?
The main purpose of the present paper is to give an answer to this question.

In Section \ref{sec:NAmn}, for each $m, n \in \mathbb{N}$, the logic $\NAmn$ is introduced by adding the axiom scheme $\Amn$ into $\N$. 
We then introduce the notion of \accessibility of $\N$-frames, and prove that $\NAmn$ is indeed sound with respect to \accessible $\N$-frames. 

The logic $\NAmn$ is expected to be complete with respect to \accessible $\N$-frames, but interestingly, it is not the case in general. 
Actually, in Section \ref{sec:NRAmn}, we prove that for $n \ge 2$, the logic $\NA{0}{n}$ is incomplete with respect to \accessible[0,n] $\N$-frames. 
Here, we pay attention to the weak variant $\RosBox$: $\dfrac{\neg \Box \varphi}{\neg \Box \Box \varphi}$ of the Rosser rule \textsc{Ros}: $\dfrac{\neg \varphi}{\neg \Box \varphi}$ which was introduced in \cite{Kurahashi2022} to analyze the properties of Rosser provability predicates.
We then introduce the logic $\NRAmn$ by adding the rule $\RosBox$ into $\NAmn$. 
In the case of $m \ge 1$ or $n \le 1$, the logics $\NAmn$ and $\NRAmn$ coincide, while in the case of $m = 0$ and $n \ge 2$, the logic $\NRA{0}{n}$ is a proper extension of $\NA{0}{n}$. 
We prove that $\NRAmn$ is also sound with respect to \accessible $\N$-frames.

In Sections \ref{sec:ffp} and \ref{sec:ffp2}, we prove that $\NRAmn$ is complete with respect to \accessible $\N$-frames. 
Moreover, we prove the finite frame property of $\NRAmn$, that is, $\NRAmn$ is characterized by the class of all finite \accessible $\N$-frames. 
Section \ref{sec:ffp} is mainly devoted to the proof of the $n \geq 1$ case, and the proof of the $n=0$ case is given in Section \ref{sec:ffp2}.
As a corollary to our proof, we obtain the decidability of $\NRAmn$. 
Finally, in Section \ref{sec:futurework}, we discuss future work.

\section{Preliminaries}\label{sec:pre}

Let $\PropVar$ denote the set of all propositional variables. 
Let $\mathscr{L}_\Box$ be the language of modal propositional logic, which consists of propositional variables $p, q, \ldots \in \PropVar$, logical constant $\bot$, propositional connectives $\neg, \lor$, and modal operator $\Box$ with the following abbreviations:
          \begin{tasks}[label=\labelitemii](2)
            \task $\top \defeq \neg \bot$
            \task $\Diamond \varphi \defeq \neg \Box \neg \varphi$
            \task $(\varphi \land \psi) \defeq \neg (\neg \varphi \lor \neg \psi)$
            \task $(\varphi \to \psi) \defeq (\neg \varphi \lor \psi)$
          \end{tasks}
Let $\MF$ denote the set of all $\mathscr{L}_\Box$-formulae.
For every $\varphi \in \MF$, let $\Sub(\varphi)$ denote the set of all the subformulae of $\varphi$.

Here we introduce the modal logic $\N$, which is originally introduced by \FMT \cite{fmt}. 
The axioms of $\N$ are propositional tautologies in the language $\mathscr{L}_\Box$ and the inference rules of $\N$ are $\dfrac{\varphi \quad \varphi \to \psi}{\psi}$ (\textsc{MP}) and $\dfrac{\varphi}{\Box \varphi}$ (\textsc{Nec}). 


Since $\N$ is non-normal, the usual Kripke semantics cannot be used with $\N$. 
\FMT introduced a Kripke-like semantics for $\N$ with some modification
on the accessibility relations.

\begin{definition}[$\N$-frames and $\N$-models {\cite[Section 3]{fmt}}] \label{def:n-model} \leavevmode
  \begin{itemize}
    \item $(W, \Rels)$ is an $\N$-\emph{frame} iff $W$ is a non-empty set and $\Rel{\varphi}$ is a
          binary relation on $W$ for every $\varphi \in \MF$.
    \item $(\mathcal{F}, V)$ is an $\N$-\emph{model} \emph{based on} $\mathcal{F}$
          iff $\mathcal{F}$ is an $\N$-frame and $V$ is a function from $W \times \PropVar$ to $\{0,1\}$.
    \item Let $M = (W, \Rels, V)$ be any $\N$-model.
          We define the \emph{satisfaction relation} $w \Vdash_M \psi$
          on $W \times \MF$ by induction on the construction of $\psi$ as follows:
          \begin{itemize}
            \item $w \Vdash_M p$ iff $V(w, p) = 1$.
            \item $w \nVdash_M \bot$.
          \item $w \Vdash_M \neg \psi$ iff $w \nVdash_M \psi$.
            \item $w \Vdash_M \psi_1 \lor \psi_2$ iff ($w \Vdash_M \psi_1$ or $w \Vdash_M \psi_2$).
            \item $w \Vdash_M \Box \psi$ iff $w' \Vdash_M \psi$ for every $w' \in W$ such that $w \Rel{\psi} w'$.
          \end{itemize}
    \item We write $(\mathcal{F}, \Vdash)$ to mean an $\N$-model $M = (\mathcal{F}, V)$
          with its satisfaction relation $\Vdash \;=\; \Vdash_M$ defined as above,
          when we do not need to talk about $V$.
    \item $\psi \in \MF$ is \emph{valid} in an $\N$-model $(W, \Rels, \Vdash)$
          if $w \Vdash \psi$ for every $w \in W$.
    \item $\psi \in \MF$ is \emph{valid} on an $\N$-frame $\mathcal{F}$
          if $\psi$ is valid in every $\N$-model $(\mathcal{F}, \Vdash)$ based on $\mathcal{F}$.
  \end{itemize}
\end{definition}

\FMT proved that $\N$ is sound, complete, and has the finite frame property with respect to the above semantics:

\begin{theorem}[{\cite[Theorems 3.6 and 4.10]{fmt}}] \label{thm:n-is-complete}
  For any $\psi \in \MF$, the following are equivalent:
  \begin{enumerate}
    \item $\N \vdash \psi$.
    \item $\psi$ is valid on every $\N$-frame.
    \item $\psi$ is valid on every finite $\N$-frame.
  \end{enumerate}
\end{theorem}

\FMT also showed that we have to look at only a finite subset of $\Rels$
when we are considering the truth of a single formula $\psi \in \MF$:

\begin{proposition}[{\cite[Theorem 4.11]{fmt}}]
  Let $\psi \in \MF$. 
  Let $M = (W, \Rels, V)$ and $M' = (W,\Rels['],V)$ be $\N$-models
  such that $\Rel{\varphi} = \Rel{\varphi}'$ for every $\Box \varphi \in \Sub(\psi)$, then for every $w \in W$, we have $w \Vdash_M \psi$ if and only if $w \Vdash_{M'} \psi$.
\end{proposition}

\begin{corollary} \label{cor:sub-preserves-validity}
  Let $\psi \in \MF$. 
  Let $\mathcal{F} = (W, \Rels)$ and $\mathcal{F}' = (W, \Rels['])$ be $\N$-frames
  such that $\Rel{\varphi} = \Rel{\varphi}'$ for every $\Box \varphi \in \Sub(\psi)$.
  Then $\psi$ is valid on $\mathcal{F}$ if and only if $\psi$ is valid on $\mathcal{F}'$.
\end{corollary}

Recently in \cite{Kurahashi2022}, the three extensions $\mathbf{NR}$, $\mathbf{N4}$, and $\mathbf{NR4}$ of $\N$ were introduced, and the finite frame property of these logics was proved. 

\begin{definition}\leavevmode
\begin{itemize}
	\item The logic $\mathbf{NR}$ is obtained from $\N$ by adding the inference rule $\dfrac{\neg \varphi}{\neg \Box \varphi}$ which is called the \emph{Rosser rule} (\textsc{Ros}). 
	\item The logics $\mathbf{N4}$ and $\textbf{NR4}$ are obtained from $\N$ and $\mathbf{NR}$, respectively, by adding the axiom scheme $\Box \varphi \to \Box \Box \varphi$. 
\end{itemize}
\end{definition}

\begin{definition}
Let $\psi \in \MF$ and $\Gamma \subseteq \MF$. 
Let $\mathcal{F} = (W, \Rels)$ be any $\N$-frame. 
\begin{itemize}
	\item $\mathcal{F}$ is \emph{$\psi$-serial} iff for every $w \in W$, there exists $x \in W$ such that $w \Rel{\psi} x$. 
	\item $\mathcal{F}$ is \emph{$\Gamma$-serial} iff $\mathcal{F}$ is $\psi$-serial for every $\Box \psi \in \Gamma$.  
	\item $\mathcal{F}$ is \emph{serial} iff $\mathcal{F}$ is $\MF$-serial. 
	\item $\mathcal{F}$ is \textit{$\psi$-transitive} iff for every $x, y, z \in W$, if $x \Rel{\Box \psi} y$ and $y \Rel{\psi} z$, then $x \Rel{\psi} z$. 
	\item $\mathcal{F}$ is said to be \textit{$\Gamma$-transitive} iff $\mathcal{F}$ is $\psi$-transitive for every $\Box \Box \psi \in \Gamma$.  
	\item $\mathcal{F}$ is called \textit{transitive} iff $\mathcal{F}$ is $\MF$-transitive. 
\end{itemize}
\end{definition}

\begin{theorem}[The finite frame property of $\mathbf{NR}$ {\cite[Theorem 3.12]{Kurahashi2022}}]\label{thm:complNR}
For any $\psi \in \MF$, the following are equivalent: 
\begin{enumerate}
	\item $\mathbf{NR} \vdash \psi$. 
	\item $\psi$ is valid on all serial $\N$-frames. 
	\item $\psi$ is valid on all finite serial $\N$-frames. 
	\item $\psi$ is valid on all finite $\Sub(\psi)$-serial $\N$-frames. 
\end{enumerate}
\end{theorem}

\begin{theorem}[The finite frame property of $\mathbf{N4}$ {\cite[Theorem 3.13]{Kurahashi2022}}]\label{thm:complNF}
For any $\psi \in \MF$, the following are equivalent: 
\begin{enumerate}
	\item $\mathbf{N4} \vdash \psi$. 
	\item $\psi$ is valid on all transitive $\N$-frames. 
	\item $\psi$ is valid on all finite transitive $\N$-frames. 
	\item $\psi$ is valid on all finite $\Sub(\psi)$-transitive $\N$-frames. 
\end{enumerate}
\end{theorem}

\begin{theorem}[The finite frame property of $\mathbf{NR4}$ {\cite[Theorem 3.14]{Kurahashi2022}}]\label{thm:complNRF}
For any $\psi \in \MF$, the following are equivalent: 
\begin{enumerate}
	\item $\mathbf{NR4} \vdash \psi$. 
	\item $\psi$ is valid on all transitive and serial $\N$-frames. 
	\item $\psi$ is valid on all finite transitive and serial $\N$-frames. 
	\item $\psi$ is valid on all finite $\Sub(\psi)$-transitive and $\Sub(\psi)$-serial $\N$-frames. 
\end{enumerate}
\end{theorem}

By using Theorems \ref{thm:complNR}, \ref{thm:complNF}, and \ref{thm:complNRF}, it was proved in \cite{Kurahashi2022} that the logics $\mathbf{NR}$, $\mathbf{N4}$, and $\mathbf{NR4}$ are exactly the provability logics of all Rosser provability predicates, all provability predicates satisfying the condition $\mathsf{D3}$, and all Rosser provability predicates satisfying $\mathsf{D3}$, respectively. 

In this paper, we focus on the logic $\mathbf{N4}$ and attempt to generalize Theorem \ref{thm:complNF}. 
We leave the analysis of logics with the rule \textsc{Ros} to future work (See Section \ref{sec:futurework}).

\section{The logic $\NAmn$}\label{sec:NAmn}

In this section, we introduce infinitely many natural extensions $\NAmn$ of $\N$, including $\mathbf{N4}$, and prove the soundness of these logics with respect to the corresponding $\N$-frames.

Let $m, n \in \mathbb{N}$.
The logic $\NAmn$ is obtained from $\N$ by adding the following axiom scheme $\Amn$:
$$\Amn :\ \Box^n \varphi \to \Box^m \varphi$$
The logic $\NA{2}{1}$ is exactly $\mathbf{N4}$. 
As a generalization of the notion of transitivity, we introduce that of \emph{\accessibility}.

\begin{definition}[\accessibility] \label{def:mn-accessibility}
  Let $\mathcal{F} = (W, \Rels)$ be any $\N$-frame.
  \begin{itemize}
    \item A path
        \[
            w_k \Rel{\Box^{k-1} \varphi} w_{k-1} \Rel{\Box^{k-2} \varphi} \cdots \Rel{\Box \varphi} w_1 \Rel{\varphi} w_0
        \]
          of elements of $W$ is called a $\varphi$-\emph{path of length k} from $w_k$ to $w_0$.

    \item Let $x, y \in W$. We write $x \Rel{\varphi}^k y$ to mean that there is a $\varphi$-path of length $k$ from $x$ to $y$.
          More formally:
          \begin{eqnarray*}
            x \Rel{\varphi}^k y
            &\defsiff&
            \left\{\begin{array}{ll}
              x = y                                                             & \tif k = 0,  \\
              \some{w \in W} x \Rel{\Box^{k-1} \varphi} w \Rel{\varphi}^{k-1} y & \tif k \ge 1.
            \end{array}
            \right. 
          \end{eqnarray*}

    \item Let $\psi \in \MF$. 
    We say $\mathcal{F}$ is $\psi$-\emph{\accessible}
          iff for any $x, y \in W$, if $x \Rel{\psi}^m y$, then $x \Rel{\psi}^n y$. 
    \item Let $\Gamma \subseteq \MF$. 
    We say $\mathcal{F}$ is $\Gamma$-\emph{\accessible} iff $\mathcal{F}$ is $\psi$-\accessible
          for every $\psi \in \MF$ such that $\Box^m \psi \in \Gamma$.

    \item We say $\mathcal{F}$ is \emph{\accessible} iff $\mathcal{F}$ is $\MF$-\accessible.
  \end{itemize}
\end{definition}

Note that \accessibility[2,1] coincides with transitivity. 
The following proposition is proved easily. 

\begin{proposition} \label{prop:value-of-box-k}
  Let $(W, \Rels, \Vdash)$ be any $\N$-model, $w \in W$, $\psi \in \MF$, and $k \ge 0$, then $w \Vdash \Box^k \psi$ if and only if $w' \Vdash \psi$ for all $w' \in W$ such that $w \Rel{\psi}^k w'$.
\end{proposition}

We are ready to prove the soundness of $\NA{m}{n}$ with respect to \accessible $\N$-frames.

\begin{theorem}[Soundness of $\NAmn$] \label{thm:soundness}
  Let $\psi \in \MF$. 
  If $\NAmn \vdash \psi$, then $\psi$ is valid on every \accessible $\N$-frame.
\end{theorem}
\begin{proof}
    Since $\N$ is sound with respect to $\N$-frames by Theorem \ref{thm:n-is-complete}, it suffices to show that
    the axiom $\Amn: \Box^n \rho \to \Box^m \rho$ is valid on every \accessible $\N$-frame.

    Let $(W, \Rels, \Vdash)$ be any $\N$-model with an \accessible $\N$-frame and $\rho \in \MF$. 
    We take any $w \in W$ and show that $w \Vdash \Box^n \rho \to \Box^m \rho$.
    Suppose that $w \Vdash \Box^n \rho$.
    Let $w' \in W$ be such that $w \Rel{\rho}^m w'$,
    then by the \accessibility, we have $w \Rel{\rho}^n w'$.
    This and $w \Vdash \Box^n \rho$ imply that $w' \Vdash \rho$ by Proposition \ref{prop:value-of-box-k}.
    Therefore, $w \Vdash \Box^m \rho$ by Proposition \ref{prop:value-of-box-k}. 
\end{proof}

\section{The rule $\RosBox$ and the logic $\NRAmn$}\label{sec:NRAmn}

We would expect the logic $\NAmn$ to be complete with respect to \accessible $\N$-frames, but let us pause here.
We pay attention to the following weak variant $\RosBox$ of the Rosser rule \textsc{Ros}: 

$$\RosBox:\ \dfrac{\neg \Box \varphi}{\neg \Box \Box \varphi}$$

For $m,n \in \mathbb{N}$, the logic $\NRAmn$ is obtained from $\NAmn$ by adding the rule $\RosBox$. 
It is easily shown that the rule $\RosBox$ is admissible in logics which is closed under the rules \textsc{Nec} and $\dfrac{\varphi \to \psi}{\Box \varphi \to \Box \psi}$ (\textsc{M}). 
Our logic $\NAmn$ is not necessarily closed under the rule \textsc{M}, but we first show that the rule $\RosBox$ is trivially admissible in many cases.

\begin{proposition} \label{prop:unprovability}
  If $m \ge 1$, then $\NAmn \nvdash \neg \Box \psi$ for all $\psi \in \MF$. 
\end{proposition}
  \begin{proof}
  Let $\mathcal{F} = (W, \Rels)$ be the $\N$-frame defined as follows: 
  \begin{itemize}
    \item $W = \{a\}$. 
    \item $\Rel{\varphi} = \emptyset$ for all $\varphi \in \MF$. 
  \end{itemize}
  Since there is no $x \in W$ such that $a \Rel{\varphi}^m x$ for every $\varphi \in \MF$, the frame $\mathcal{F}$ is trivially \accessible. 
  For every $\psi \in \MF$, since there is no $x \in W$ such that $a \Rel{\psi} x$, we have that $\Box \psi$ is valid on $\mathcal{F}$. 
  By Theorem \ref{thm:soundness}, we conclude that $\NAmn \nvdash \neg \Box \psi$. 
  \end{proof}

\begin{proposition}\label{prop:admissibility}
  Suppose $m \ge 1$ or $n \leq 1$. 
  Then, the rule $\RosBox$ is admissible in $\NAmn$. 
  Consequently, we have $\NRAmn = \NAmn$. 
\end{proposition}
\begin{proof}
    We distinguish the following three cases: 
    \begin{itemize}
        \item Case 1: $m \ge 1$. \\
        The admissibility of $\RosBox$ in $\NAmn$ immediately follows from Proposition \ref{prop:unprovability}. 
        \item Case 2: $m = 0$ and $n = 0$. \\
        Since $\NA{0}{0} = \N = \NA{1}{1}$, the admissibility of $\RosBox$ in $\NA{0}{0}$ is shown in Case 1. 
        \item Case 3: $m = 0$ and $n = 1$. \\
        Since $\NA{0}{1} \vdash \neg \Box \varphi \to \neg \Box \Box \varphi$, the rule $\RosBox$ is also trivially admissible in $\NA{0}{1}$. \qedhere
    \end{itemize}
  \end{proof}

Thus, the remaining case is $m = 0$ and $n \ge 2$. 
The following proposition shows that in this case, $\RosBox$ is an important rule with respect to \accessible[0, n] $\N$-frames.

\begin{proposition} \label{prop:soundness-rosbox}
  Let $\psi \in \MF$ and $n \ge 2$. 
  \begin{enumerate}
    \item If $\psi$ is not of the form $\Box^{n-1} \varphi$ for every $\varphi \in \MF$, then there is a \accessible[0,n] $\N$-frame $M$ in which $\neg \Box \psi$ is not valid.
    \item If $\neg \Box^n \psi$ is valid in a \accessible[0,n] $\N$-model $M = (W, \Rels, V)$, then so is $\neg \Box^{n+1} \psi$.
  \end{enumerate}
\end{proposition}
\begin{proof}
   1. Assume that $\psi$ is not of the form $\Box^{n-1} \varphi$ for every $\varphi \in \MF$. 
    We define the $\N$-model $M = (W, \Rels, V)$ as follows:
    \begin{itemize}
      \item $W = \set{ a, b }$;
      \item For $w_1, w_2 \in W$ and $\varphi \in \MF$, $w_1 \Rel{\varphi} w_2$ iff ($w_1 \neq a$ or $\varphi \neq \psi$);
      \item $V$ is arbitrary.
    \end{itemize}

The model $M$ is visualized in Figure \ref{fig0}.    

\begin{figure}[th]
\centering
\begin{tikzpicture}
\node [draw, circle] (a) at (0,0) {$a$};
\node [draw, circle] (b) at (0,2) {$b$};

\draw [thick, ->] (b) to [out=225,in=135] (a);
\draw [thick, ->] (a) to [out=45,in=315] (b);
\draw [thick, ->] (0.2,2.2) arc (-45:225:0.3);

\draw (0.1, 2.8) node[right] {$\Rel{\varphi}$};
\draw (0.5, 1) node[right] {$\Rel{\varphi}$ ($\varphi \ne \psi$)};
\draw (-0.5, 1) node[left] {$\Rel{\varphi}$};

\end{tikzpicture}
\caption{The model $M$}
    \label{fig0}
\end{figure}
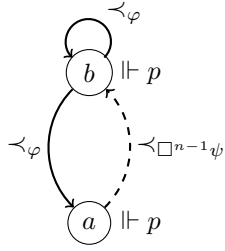    
    
    Here, $a \Vdash_{M} \Box \psi$ since there is no $w \in W$ such that $a \Rel{\psi} w$.
    So it suffices to show that $(W, \Rels)$ is \accessible[0,n].
    Let $\varphi \in \MF$. 
    It is easy to see that $b \Rel{\varphi}^n b$. 
    Also, we have $a \Rel{\Box^{n-1} \varphi} b$ because $\psi \neq \Box^{n-1} \varphi$. 
    Thus, 
    $$
              a \Rel{\Box^{n-1} \varphi} b \Rel{\Box^{n-2} \varphi} b \Rel{\Box^{n-3} \varphi} \cdots \Rel{\Box \varphi} b \Rel{\varphi} a.
    $$
    Therefore, we get $a \Rel{\varphi}^{n} a$. 
    We have proved that $(W, \Rels)$ is \accessible[0,n].

   2. We prove the contrapositive.
    Assume that $w^* \Vdash_M \Box^{n+1} \psi$ for some $w^* \in W$.
    Here, \accessibility[0,n] of $(W, \Rels)$ implies that
    $w^* \Rel{\varphi}^n w^*$ for any $\varphi \in \MF$.
    Then, $w^* \Rel{\Box \psi}^n w^*$ implies that
    there are $w_1, w_2, \ldots, w_{n-1}\in W$ such that:
    $$
      w^* \Rel{\Box^n \psi} w_{n-1} \Rel{\Box^{n-1} \psi} \cdots \Rel{\Box^2 \psi} w_1 \Rel{\Box \psi} w^* .
    $$
    Here, $w^* \Vdash_M \Box^{n+1} \psi$ and $w^* \Rel{\Box^n \psi} w_{n-1}$ imply that $w_{n-1} \Vdash_M \Box^n \psi$. 
    Therefore, $\neg \Box^n \psi$ is also not valid in $M$.
 \end{proof}

We obtain the following refinement of Theorem \ref{thm:soundness}. 

\begin{theorem}[Soundness of $\NRAmn$] \label{thm:soundness2}
  Let $\psi \in \MF$. 
  If $\NRAmn \vdash \psi$, then $\psi$ is valid on every \accessible $\N$-frame.
\end{theorem}
\begin{proof}
By Theorem \ref{thm:soundness} and Proposition \ref{prop:admissibility}, it suffices to show the theorem for the case $m = 0$ and $n \ge 2$. 
We prove the theorem by induction on the length of proofs in $\NRA{0}{n}$. 
By the proof of Theorem \ref{thm:soundness}, it suffices to show that if $\NRA{0}{n} \vdash \neg \Box \psi$ and $\neg \Box \psi$ is valid on all \accessible[0,n] $\N$-frames, then $\neg \Box \Box \psi$ is also valid on all \accessible[0,n] $\N$-frames. 
This is an immediate consequence of Proposition \ref{prop:soundness-rosbox}. 
\end{proof}

We have $\NA{0}{n} \vdash \neg \Box^n \bot$ because $\Box^n \bot \to \bot$ is an instance of $\Amn[0,n]$. 
By Theorem \ref{thm:soundness}, the formula $\neg \Box^n \bot$ is valid on all \accessible[0,n] $\N$-frames. 
Then, by Proposition \ref{prop:soundness-rosbox}, we obtain that the formula $\neg \Box^{n+1} \bot$ is also valid on all \accessible[0,n] $\N$-frames. 
On the other hand, the following proposition shows that $\neg \Box^{n+1} \bot$ is not provable in $\NA{0}{n}$ for $n \ge 2$. 

\begin{proposition}
For $n \ge 2$, $\NA{0}{n} \nvdash \neg \Box^{n+1} \bot$. 
Consequently, $\NA{0}{n} \subsetneq \NRA{0}{n}$. 
\end{proposition}
\begin{proof}
    It suffices to construct an $\N$-model in which $\Amn[0,n]: \Box^n \psi \to \psi$ is valid for every $\psi \in \MF$ but $\neg \Box^{n+1} \bot$ is not.
    Let $M = (\set{a, b}, \Rels, V)$ be the $\N$-model defined as follows: 
    \begin{itemize}
        \item $V(a, p) = V(b, p) = 1$ for any $p \in \PropVar$;
        \item $b \Rel{\varphi} w$ for any $\varphi \in \MF$ and any $w \in \set{a,b}$;
        \item $a \not \Rel{\varphi} a$ for any $\varphi \in \MF$;
        \item If $\varphi$ is not of the form $\Box^{n-1} \psi$, then $a \not \Rel{\varphi} b$; 
        \item $a \Rel{\Box^{n-1} \varphi} b$ is defined by induction on $\varphi$:
        \begin{itemize}
            \item $a \Rel{\Box^{n-1} \bot} b$.
            \item $a \not \Rel{\Box^{n-1} p} b$ for any $p \in \PropVar$.
            \item $a \Rel{\Box^{n-1} \neg \psi} b$ iff $a \not \Rel{\Box^{n-1} \psi}b$.
            \item $a \Rel{\Box^{n-1} (\psi_1 \lor \psi_2)} b$ iff $a \Rel{\Box^{n-1} \psi_1} b$ and $a \Rel{\Box^{n-1} \psi_2} b$.
            \item $a \Rel{\Box^{n-1} \Box \psi} b$ iff $a \Rel{\psi} b$.
        \end{itemize}
    \end{itemize}
    Note that the definition of $a \Rel{\varphi} b$ makes sense because $n \ge 2$. 
    The model $M$ is visualized in Figure \ref{fig1}. 

\begin{figure}[th]
\centering
\begin{tikzpicture}
\node [draw, circle] (a) at (0,0) {$a$};
\node [draw, circle] (b) at (0,2) {$b$};

\draw [thick, ->] (b) to [out=225,in=135] (a);
\draw [thick, ->, dashed] (a) to [out=45,in=315] (b);
\draw [thick, ->] (0.2,2.2) arc (-45:225:0.3);

\draw (0.3, 0) node[right] {$\Vdash p$};
\draw (0.3, 2) node[right] {$\Vdash p$};
\draw (0.1, 2.8) node[right] {$\Rel{\varphi}$};
\draw (0.5, 1) node[right] {$\Rel{\Box^{n-1} \psi}$};
\draw (-0.5, 1) node[left] {$\Rel{\varphi}$};

\end{tikzpicture}
\caption{The model $M$}
    \label{fig1}
\end{figure}

The dashed line in the figure indicates that $a \Rel{\Box^{n-1} \psi} b$ holds not necessarily for all $\psi$. 
The following claim concerns the condition for $a \Rel{\Box^{n-1} \psi} b$ to be held.

    \begin{claim*}
        For any $\psi \in \MF$, $a \Rel{\Box^{n-1} \psi} b$ if and only if $a \nVdash \psi$.
        \begin{subproof}[Proof of the claim]
            We prove the claim by induction on the construction of $\psi$.
            \begin{itemize}
                \item $a \Rel{\Box^{n-1} \bot} b$ and $a \nVdash \bot$.
                \item $a \not \Rel{\Box^{n-1} p} b$ and $a \Vdash p$ for any $p \in \PropVar$.
                \item $\begin{aligned}[t]
                a \Rel{\Box^{n-1} \neg \psi'} b
                    &\iff a \not \Rel{\Box^{n-1} \psi'} b \\
                    &\iff a \Vdash \psi' & \textrm{(by I.H.)} \\
                    &\iff a \nVdash \neg \psi'.
                \end{aligned}$
                \item $\begin{aligned}[t]
                a \Rel{\Box^{n-1} (\psi_1 \lor \psi_2)} b
                    &\iff a \Rel{\Box^{n-1} \psi_1} b \And a \Rel{\Box^{n-1} \psi_2} b \\
                    &\iff a \nVdash \psi_1 \And a \nVdash \psi_2 & \textrm{(by I.H.)} \\
                    &\iff a \nVdash \psi_1 \lor \psi_2.
                \end{aligned}$
                \item We shall prove $a \Rel{\Box^{n-1} \Box \psi'} b \siff a \nVdash \Box \psi'$ in both directions.
                \begin{itemize}
                    \item[$(\Rightarrow)$]
                    Assume that $a \Rel{\Box^{n-1} \Box \psi'} b$, then $a \Rel{\psi'} b$ by the definition of $\Rel{\Box^{n-1} \Box \psi'}$.
                    Here, $a \Rel{\psi'} b$ implies that $\psi'$ should be of the form $\Box^{n-1} \psi''$. Now that $a \Rel{\Box^{n-1} \psi''} b$, then $a \nVdash \psi''$ by the induction hypothesis. 
                    Recall that $b \Rel{\varphi} w$ for any $\varphi \in \MF$ and any $w \in \set{a, b}$, then we have $a \Rel{\psi''}^n a$ because the following path exists:
                    $$
                    a \Rel{\Box^{n-1} \psi''} b \Rel{\Box^{n-2} \psi''} b \Rel{\Box^{n-2} \psi''} \cdots \Rel{\Box \psi''} b \Rel{\psi''} a.
                    $$
                    This implies that $a \nVdash \Box^n \psi''$, that is, $a \nVdash \Box \psi'$.
                    \item[$(\Leftarrow)$]
                    Assume that $a \nVdash \Box \psi'$, then $a \Rel{\psi'} b$ (and $b \nVdash \psi'$).
                    This implies that $a \Rel{\Box^{n-1} \Box \psi'} b$. \qedhere
                \end{itemize}
            \end{itemize}
        \end{subproof}
    \end{claim*}
    Here, $a \not \Rel{\bot} b$ implies $a \not \Rel{\Box^n \bot} b$ by the definition of $\Rel{\Box^n \bot}$. 
    Since there is no $w \in W$ such that $a \Rel{\Box^n \bot} w$, we obtain $a \Vdash \Box^{n+1} \bot$. 
    Therefore, $\neg \Box^{n+1} \bot$ is not valid in $M$.

    Now it suffices to show that $w \Vdash \Box^n \psi \to \psi$ for every $\psi \in \MF$ and $w \in \set{a, b}$.
    Take any $\psi \in \MF$. 
    We shall first assume that $a \nVdash \psi$ and show that $a \nVdash \Box^n \psi$. Here, $a \nVdash \psi$ implies that $a \Rel{\Box^{n-1} \psi} b$ by the claim, then $a \Rel{\Box^{2n-1} \psi} b$ by the definition of $a \Rel{\varphi} b$, so $a \nVdash \Box^n \psi$ by the claim.
    Now we shall show that $b \Vdash \Box^n \psi \to \psi$. Since $b \Rel{\varphi} b$ holds for any $\varphi \in \MF$, we have $b \Rel{\psi}^n b$. 
    Hence, $b \Vdash \Box^n \psi$ trivially implies $b \Vdash \psi$, and so $b \Vdash \Box^n \psi \to \psi$ holds. 
\end{proof}

\begin{corollary}
For $n \ge 2$, the logic $\NA{0}{n}$ is incomplete with respect to \accessible[0,n] $\N$-frames. 
\end{corollary}

\section{The finite frame property of $\NRAmn$}\label{sec:ffp}

In the following two sections, we prove the completeness and the finite frame property of the logics $\NRAmn$.
More precisely, we prove the following main theorem of the present paper. 

\begin{theorem}\label{thm:ffp}
  For any $\psi \in \MF$, the following are equivalent:
  \begin{enumerate}
    \item $\NRAmn \vdash \psi$.
    \item $\psi$ is valid on every \accessible $\N$-frame.
    \item $\psi$ is valid on every finite \accessible $\N$-frame.
    \item $\psi$ is valid on every finite $\Sub(\psi)$-\accessible $\N$-frame.
  \end{enumerate}
\end{theorem}

Here, (1 $\Rightarrow$ 2) immediately follows from Theorem \ref{thm:soundness2}, and (2 $\Rightarrow$ 3) is obvious.
We prove the implication (3 $\Rightarrow$ 4) in the following Lemma \ref{lem:3to4}. 
The proof of the implication $(4 \Rightarrow 1)$ is lengthy, and so we divide it into two parts. 
The first part of the proof is the case of $n \geq 1$ and is proved in this section. 
We prove the second part, which is the case of $n = 0$, in the next section. 
We first prove the implication (3 $\Rightarrow$ 4). 

\begin{lemma}[(3 $\Rightarrow$ 4) of Theorem \ref{thm:ffp}] \label{lem:3to4}
  For any $\psi \in \MF$, if $\psi$ is valid on all finite \accessible $\N$-frames, then
  $\psi$ is also valid on all finite $\Sub(\psi)$-\accessible $\N$-frames.
\end{lemma}
\begin{proof}
    Assume that $\psi$ is valid on all finite \accessible $\N$-frames.
    Let $\mathcal{F} = (W, \Rels)$ be any finite $\Sub(\psi)$-\accessible $\N$-frame.
    We would like to prove that $\psi$ is valid on $\mathcal{F}$. 
    We define $\mathcal{F^\ast} = (W, \Rels[^\ast])$ as follows:
    \begin{enumerate}[label={(\roman*)}]
      \item \label{3to4:relstar1}
            if $\Box \varphi \in \Sub(\psi)$, then $x \Rel{\varphi}^\ast y$ iff $x \Rel{\varphi} y$;
      \item \label{3to4:relstar2}
            if $\varphi \in \Sub(\psi)$, $\Box \varphi \notin \Sub(\psi)$, and $n > m$,
            then $x \Rel{\varphi}^* y$;
      \item \label{3to4:relstar3}
            if $\varphi \in \Sub(\psi)$, $\Box \varphi \notin \Sub(\psi)$, and $m > n$,
            then $x \not \Rel{\varphi}^* y$;
      \item \label{3to4:relstar4} otherwise, $x \Rel{\varphi}^\ast y$ iff $x = y$.
    \end{enumerate}
    We shall show that $\mathcal{F^\ast}$ is \accessible; then $\psi$ is valid on $\mathcal{F^\ast}$ by the assumption. 
    Since $\Rel{\varphi} = \Rel{\varphi}^\ast$ for every $\Box \varphi \in \Sub(\psi)$, by Proposition \ref{cor:sub-preserves-validity}, we conclude that $\psi$ is also valid on $\mathcal{F}$.

    Take any $\rho \in \MF$ and any $x_0, x_1, \ldots, x_{m-1} \in W$ such that:
    \begin{equation}\label{3to4:path}
      x_m \Rel{\Box^{m-1} \rho}^\ast x_{m-1} \Rel{\Box^{m-2} \rho}^\ast \cdots \Rel{\Box \rho}^\ast x_1 \Rel{\rho}^\ast x_0 .
    \end{equation}
    Then we would like to show that $x_m \Rel{\rho}^{*n} x_0$.

    If $\rho \notin \Sub(\psi)$, then (\ref{3to4:path}) and \ref{3to4:relstar4} imply that $x_m = x_{m-1} = \ \cdots\ = x_0$.
    Here, \ref{3to4:relstar4} also implies that $x_m \Rel{\Box^i \rho}^* x_m$ for any $i \ge 0$.
    Then it follows that $x_m \Rel{\rho}^{*n} x_0$ because the following path exists:
    $$
      x_m \Rel{\Box^{n-1}\rho}^* x_m \Rel{\Box^{n-2}\rho}^* \cdots \Rel{\Box \rho}^* x_m \Rel{\rho}^* x_m = x_0 .
    $$
    So we may assume that $\rho \in \Sub(\psi)$. 
    Let $k = \max \set{ i; \Box^i \rho \in \Sub(\psi) }$.
    We distinguish the following two cases. 

    \paragraph*{Case 1:} $k \ge m$. \\  
    Since $\Box^m \rho \in \Sub(\psi)$, (\ref{3to4:path}) and \ref{3to4:relstar1} imply that:
    $$
      x_m \Rel{\Box^{m-1} \rho} x_{m-1} \Rel{\Box^{m-2} \rho} \cdots \Rel{\Box \rho} x_1 \Rel{\rho} x_0 .
    $$
    Then by the \accessibility of $\mathcal{F}$, there are $y_0, y_1, \ldots, y_{n-1} \in W$ such that:
    \begin{equation}\label{3to4:case1}
      x_m \Rel{\Box^{n-1} \rho} y_{n-1} \Rel{\Box^{n-2} \rho} \cdots \Rel{\Box \rho} y_1 \Rel{\rho} y_0 = x_0 .
    \end{equation}
    If $k \ge n$, we have $\Box^n \rho \in \Sub(\psi)$, and then (\ref{3to4:case1}) and \ref{3to4:relstar1} imply that $x_m \Rel{\rho}^{*n} x_0$.
    So we may assume that $n > k$. 
    \begin{itemize}
      \item Since $\Box^{k+1} \rho \notin \Sub(\psi)$, \ref{3to4:relstar4} implies that
            \[
                x_m \Rel{\Box^{n-1}\rho}^* x_m \Rel{\Box^{n-2}\rho}^* \cdots \Rel{\Box^{k+1}\rho}^* x_m.
            \]
      \item Since $\Box^k \rho \in \Sub(\psi)$, (\ref{3to4:case1}) and \ref{3to4:relstar1} imply that
            \[
             y_k \Rel{\Box^{k-1}\rho}^* y_{k-1} \Rel{\Box^{k-2}\rho}^* \cdots \Rel{\Box \rho}^* y_1 \Rel{\rho}^* x_0.
            \]
      \item Since $\Box^{k+1} \rho \notin \Sub(\psi)$, $\Box^k \rho \in \Sub(\psi)$, and $n > (k \ge\,)\, m$,
            \ref{3to4:relstar2} implies that $x_m \Rel{\Box^k \rho}^* y_k$.
    \end{itemize}
    Therefore, it follows that $x_m \Rel{\rho}^{*n} x_0$ by combining these three paths.

    \paragraph*{Case 2:} $m > k$. \\   
    If $m > n$, then (iii) implies that $x \not \Rel{\Box^k \rho}^* y$ for any $x, y \in W$.
    This contradicts $x_{k+1} \Rel{\Box^k \rho}^* x_k$ from (\ref{3to4:path}). 
    Therefore, we have $n \ge m$.
    Now that $n \ge m > k$, then \ref{3to4:relstar4} implies that $x_m \Rel{\Box^i \rho}^* x_m$ for any $i \ge m$.
    Therefore, it follows that $x_m \Rel{\rho}^{*n} x_0$ because the following path exists:
    $$
      x_m \Rel{\Box^{n-1}\rho}^* x_m \Rel{\Box^{n-2}\rho}^* \cdots \Rel{\Box^m \rho}^* x_m \Rel{\Box^{m-1} \rho}^* x_{m-1} \Rel{\Box^{m-2} \rho}^* \cdots \Rel{\Box \rho}^* x_1 \Rel{\rho}^* x_0.
    $$
  \end{proof}

Next, we prove the implication (4 $\Rightarrow$ 1). 
To do that, given any $\psi \in \MF$ such that $\NRAmn \nvdash \psi$, we need to construct a finite \accessible model that falsifies $\psi$.

\begin{definition} \label{def:canonical}
  Let $\psi \in \MF$ and $X \subseteq \MF$.

  \begin{enumerate}[label=(\arabic*)]
    \item \label{def:canonical/negg}
          We let $\negg \psi \defeq \psi'$ if $\psi$ is of the form $\neg \psi'$, and $\negg \psi \defeq \neg \psi$ otherwise.
    \item \label{def:canonical/nsub}
         $\NSub{\psi} \defeq \Sub(\psi) \union \set{ \negg \rho; \rho \in \Sub(\psi) }$.
  \end{enumerate}

  Note that $\NSub{\psi}$ is finite.

  \begin{enumerate}[resume, label=(\arabic*)]
    \item \label{def:canonical/consistent} We say $X$ is \emph{$\NRAmn$-consistent} iff $\NRAmn \nvdash \neg \bigwedge X$.
    \item \label{def:canonical/maximal} We say $X$ is \emph{$\psi$-maximal} iff for every $\rho \in \NSub{\psi}$ either $\rho \in X$ or $\negg \rho \in X$.
    \item \label{def:canonical/model} We define the $\N$-model $\CanonicalModel{m,n}(\psi) = (W^*, \Rels[^*], V^*)$ as follows:
          \begin{itemize}
            \item $W^* \defeq \set{ X \subseteq \NSub{\psi}; X \textrm{ is } \psi\textrm{-maximal } \NRAmn\textrm{-consistent} }$;
            \item For every $X, Y \in W$, $X \Rel{\varphi}^* Y \defiff \Box \varphi \notin X \tor \varphi \in Y$;
            \item For every $X \in W$ and every $p \in \PropVar$, $V^*(X, p) = 1 \defiff p \in X$.
          \end{itemize}
  \end{enumerate}
\end{definition}


As usual, it is shown that every $\NRAmn$-consistent subset of $\NSub{\psi}$ can be extended to a $\psi$-maximal $\NRAmn$-consistent set. 
See \cite{cresswell2012new} for the usual arguments proving completeness and finite frame property of normal modal logics. 
The following lemma is proved in the same way as in \cite[Claim 3.17]{Kurahashi2022}. 

\begin{lemma}[Truth Lemma] \label{lem:truth-lemma}
  Let $\psi \in \MF$ and $(W^*, \Rels[^*], V^*) = \CanonicalModel{m,n}(\psi)$.
  Then for any $\rho\in \NSub{\psi}$ and any $X \in W^*$, $X \VdashCan{\psi} \rho$ iff $\rho \in X$.
  \end{lemma}

Now we are ready to prove the implication (4 $\Rightarrow$ 1) for $n \ge 1$.

\begin{lemma}[(4 $\Rightarrow$ 1) of Theorem \ref{thm:ffp} ($n \geq 1$)] \label{lem:4to1}
  Let $n \geq 1$. 
  For any $\psi \in \MF$, if $\psi$ is valid on all finite $\Sub(\psi)$-\accessible $\N$-frames, then $\NRAmn \vdash \psi$.
\end{lemma}
  \begin{proof}
    We prove the contrapositive. 
    Assume that $\NRAmn \nvdash \psi$,
    then we would like to construct a finite $\Sub(\psi)$-\accessible $\N$-frame on which $\psi$ is not valid.

    Let $(W^*, \Rels[^*], \Vdash^*) = \CanonicalModel{m,n}(\psi)$.
    Then $\NRAmn \nvdash \psi$ implies that $\set{\negg \psi}$ is $\NRAmn$-consistent.
    So, we find some $X_\psi \in W^*$ such that $\negg \psi \in X_\psi$.
    Here, $\psi \notin X_\psi$.
    Therefore by Truth Lemma, $X_\psi \nVdash^* \psi$, which implies that $\psi$ is not valid on $\mathcal{F}^* = (W^*, \Rels[^*])$.

    Now it suffices to show that $\mathcal{F}^*$ is $\Sub(\psi)$-\accessible.
    Take any \newline $\Box^m \rho \in \Sub(\psi)$ and any $x_0, x_1, \ldots, x_{m-1}, x_m \in W^*$ such that:
    \begin{equation}\label{4to1:path}
      x_m \Rel{\Box^{m-1} \rho}^\ast x_{m-1} \Rel{\Box^{m-2} \rho}^\ast \cdots \Rel{\Box \rho}^\ast x_1 \Rel{\rho}^\ast x_0 .
    \end{equation}
    We would like to show that $x_m \Rel{\rho}^{*n} x_0$.
    We distinguish the following two cases. 

    \paragraph*{Case 1:} $\Box^n \rho \in x_m$. \\
    Since $\NRAmn \vdash \Box^n \rho \to \Box^m \rho$, it is shown that $\Box^m \rho \in x_m$, 
    and so $x_m \Vdash^* \Box^m \rho$ by Truth Lemma.
    This and (\ref{4to1:path}) imply that for every $i \le m$,
    $x_i \Vdash^* \Box^i \rho$,
    and then $\Box^i \rho \in x_i$ by Truth Lemma. 
    Hence for every $i \le m$ and any $y \in W^*$, we obtain
    \begin{equation}\label{4to1:case1}
        y \Rel{\Box^i \rho}^* x_i
    \end{equation}
    by the definition of $\Rel{\Box^i \rho}^*$. 
    We distinguish the following three cases. 
    
   \subparagraph*{Case 1.1:} $m \ge n$. \\
   We have $x_m \Rel{\Box^{n-1}\rho}^* x_{n-1}$ by (\ref{4to1:case1}). 
   By (\ref{4to1:path}), we obtain that $x_m \Rel{\rho}^{*n} x_0$ because the following path exists:
    $$x_m \Rel{\Box^{n-1}\rho}^* x_{n-1} \Rel{\Box^{n-2}\rho}^*\cdots \Rel{\Box \rho}^* x_1 \Rel{\rho}^* x_0 .$$

   \subparagraph*{Case 1.2:} $m = n - 1$. \\
    By (\ref{4to1:case1}), we have $x_m \Rel{\Box^{m}\rho}^* x_m$. 
    Then (\ref{4to1:path}) implies that $x_m \Rel{\rho}^{*n} x_0$ since the following path exists:
    $$x_m \Rel{\Box^{n-1}\rho}^* x_m \Rel{\Box^{m-1}\rho}^* \cdots \Rel{\Box \rho}^* x_1 \Rel{\rho}^* x_0 .$$

   \subparagraph*{Case 1.3:} $n > m + 1$. \\
    Suppose, towards a contradiction, that $\{\Box^i \rho\}$ is $\NRAmn$-inconsistent for some $i$ with $m+ 1 \leq i < n$. 
    Then $\NRAmn \vdash \neg \Box^i \rho$, and thus $\NRAmn \vdash \neg \Box^n \rho$ by using the rule $\RosBox$ for $n-i$ times.
    This contradicts the fact $\Box^n \rho \in x_m$.

    We have proved that for each $i$ with $m+1 \leq i < n$, the set $\{\Box^i \rho\}$ is $\NRAmn$-consistent. 
    Thus, we obtain $y_{m+1}, \ldots, y_{n-1} \in W^*$ such that $\Box^i \rho \in y_i$ for each $i$. 
    By the definition of $\Rel{\Box^i \rho}^\ast$, we have 
    \begin{equation}\label{4to1:path1}
      x_m \Rel{\Box^{n-1}\rho}^* y_{n-1} \Rel{\Box^{n-2}\rho}^* \cdots \Rel{\Box^{m+1} \rho}^* y_{m+1} .
    \end{equation}
    Then since $y_{m+1} \Rel{\Box^m \rho}^* x_m$ by (\ref{4to1:case1}), it follows that $x_m \Rel{\rho}^{*n} x_0$ by connecting the paths (\ref{4to1:path1}) and (\ref{4to1:path}).

    \paragraph*{Case 2:} $\Box^n \rho \notin x_m$. \\
    By the definition of $\Rel{\Box^{n-1} \rho}^*$, for any $y \in W^*$, 
    \begin{equation}\label{4to1:case2}
         x_m \Rel{\Box^{n-1}\rho}^* y.
    \end{equation}
    We distinguish the following three cases. 

   \subparagraph*{Case 2.1:} $m \ge n$. \\
    (\ref{4to1:case2}) and (\ref{4to1:path}) imply that:
    $$x_m \Rel{\Box^{n-1}\rho}^* x_{n-1} \Rel{\Box^{n-2}\rho}^* \cdots \Rel{\Box \rho}^* x_1 \Rel{\rho}^* x_0 .$$

   \subparagraph*{Case 2.2:} $m = n - 1$. \\
    (\ref{4to1:case2}) and (\ref{4to1:path}) imply that:
    $$x_m \Rel{\Box^{n-1}\rho}^* x_m \Rel{\Box^{m-1}\rho}^* \cdots \Rel{\Box \rho}^* x_1 \Rel{\rho}^* x_0. $$

   \subparagraph*{Case 2.3:} $n > m + 1$. \\
    We shall find some $y_{m+1}, y_{m+2}, \ldots, y_{n-1} \in W^*$ such that:
    \begin{equation}\label{4to1:path2}
      y_{n-1} \Rel{\Box^{n-2}\rho}^* \cdots \Rel{\Box^{m+2} \rho}^* y_{m+2} \Rel{\Box^{m+1} \rho}^* y_{m+1} \Rel{\Box^m \rho}^* x_m .
    \end{equation}
    Then it follows that $x_m \Rel{\rho}^{*n} x_0$ by connecting the two paths (\ref{4to1:path2}) and (\ref{4to1:path}) with $x_m \Rel{\Box^{n-1}\rho}^* y_{n-1}$ by (\ref{4to1:case2}).

    Let $k = \max \set{ j; \Box^j \rho \in \NSub{\psi} }$. 
    Here, $k = \max \set { j; \Box^j \rho \in \Sub(\psi) } \ge m$. 
    We distinguish the following two cases. 
   
    \subparagraph*{Case 2.3.1:} $k \ge n$. \\
    Suppose, towards a contradiction, that the set $\set{ \neg \Box^i \rho}$ is $\NRAmn$-inconsistent for some $i$ with $m+1 \le i < n$. 
    Then, $\NRAmn \vdash \Box^i \rho$. 
    By applying the rule \textsc{Nec} for $n-i$ times, we obtain $\NRAmn \vdash \Box^n \rho$. 
    Here, $k \ge n$ implies that $\Box^n \rho \in \NSub{\psi}$.
    This contradicts the fact $\Box^n \rho \notin x_m$.
 
    We have shown that $\set{ \neg \Box^i \rho}$ is $\NRAmn$-consistent for any $i$ with $m+1 \le i < n$, then there are $y_{m+1}, y_{m+2}, \ldots, y_{n-1} \in W^*$ such that $\neg \Box^i \rho \in y_i$ for each $i$.
    Since $\Box^i \rho \notin y_i$, we have $y_i \Rel{\Box^{i-1} \rho}^* y$ for any $y \in W^*$.
    Therefore, the path (\ref{4to1:path2}) exists.

    \subparagraph*{Case 2.3.2:} $n > k \ge m$.\\
    Here, for every $i > k$, we have $\Box^i \rho \notin \NSub{\psi}$, so $\Box^i \rho \notin x_m$.
    This implies that $x_m \Rel{\Box^{i-1} \rho}^* y$ for any $y \in W^*$. 
    So we have
    $$
     x_m \Rel{\Box^{n-2}\rho}^* x_m \Rel{\Box^{n-3}\rho}^* \cdots \Rel{\Box^{k+1}\rho}^* x_m \Rel{\Box^k \rho}^* y
    $$
    for any $y \in W^*$. 
    If $k = m$, then the path (\ref{4to1:path2}) exists by letting $y = x_m$. 
    So, we may assume $k \ge m+1$. 
    It suffices to show that the rest of the path (\ref{4to1:path2}) exists, 
    that is, it suffices to show that there are $y_{m+1}, y_{m+2}, \ldots, y_k \in W^*$ such that:
    \begin{equation}\label{4to1:path3}
      y_k \Rel{\Box^{k-1} \rho}^* \cdots \Rel{\Box^{m+2} \rho}^* y_{m+2} \Rel{\Box^{m+1} \rho}^* y_{m+1} \Rel{\Box^m \rho}^* x_m .
    \end{equation}
    We distinguish the following two cases. 

    \subparagraph*{Case 2.3.2.1:} $\Box^m \rho \notin x_m$.\\
    Suppose, towards a contradiction, that the set $\set{ \neg \Box^i \rho }$ is $\NRAmn$-inconsistent for some $i$ with $m+1 \le i \le k$.
    We have $\NRAmn \vdash \Box^i \rho$. 
    By applying \textsc{Nec} for $n-i$ times, we obtain $\NRAmn \vdash \Box^n \rho$. 
    Since $\NRAmn \vdash \Box^n \rho \to \Box^m \rho$, we have $\NAmn \vdash \Box^m \rho$. 
    This contradicts $\Box^m \rho \notin x_m$.
    
    We have proved that the set $\set{ \neg \Box^i \rho }$ is $\NRAmn$-consistent for any $i$ with $m+1 \le i \le k$, then there are $y_{m+1}, y_{m+2}, \ldots, y_k \in W^*$ such that $\neg \Box^i \rho \in y_i$.
    Since $\Box^i \rho \notin y_i$, we have $y_i \Rel{\Box^{i-1} \rho}^* y$ for any $y \in W^*$.
    Therefore, the path (\ref{4to1:path3}) exists.

    \subparagraph*{Case 2.3.2.2:} $\Box^m \rho \in x_m$.\\
    We have that $y \Rel{\Box^m \rho}^* x_m$ for any $y \in W^*$. 
    We define the number $h$ as follows: 
    $$
      h = \min \left(\set{ j ; j \ge m+1 \tand \set{ \Box^j \rho } \text{ is }\NRAmn\text{-inconsistent} } \cup \{k+1\}\right).
    $$

    For every $i$ with $m+1 \le i \le h-1$, the set $\set{ \Box^i \rho }$ is $\NRAmn$-consistent,
    and thus there are $y_{m+1}, \ldots, y_{h-1} \in W^*$ such that $\Box^i \rho \in y_i$ for every $i$ with $m+1 \le i \le h-1$.
    Hence by the definition of $\Rel{\Box^{i-1}\rho}^*$, we obtain
    \begin{equation}\label{4to1:path4}
      y_{h-1} \Rel{\Box^{h-2}\rho}^* y_{h-2} \Rel{\Box^{h-3}\rho}^* \cdots \Rel{\Box^{m+1}\rho}^* y_{m+1} \Rel{\Box^m \rho}^* x_m.
    \end{equation}
    If $h = k+1$, we obtain that the path (\ref{4to1:path3}) exists.

    We may assume that $k \ge h$. 
    Then it follows that $\NRAmn \vdash \neg \Box^h \rho$.
    For every $i \ge h$, we obtain $\NRAmn \vdash \neg \Box^i \rho$ by applying $\RosBox$ for $i-h$ times. 
    Then we have $\Box^i \rho \notin x_m$,
    and thus
    $$
      x_m \Rel{\Box^{k-1} \rho}^* x_m \Rel{\Box^{k-2} \rho}^* \cdots \Rel{\Box^h \rho}^* x_m \Rel{\Box^{h-1} \rho}^* y_{h-1}. 
    $$
    By combining this with the path (\ref{4to1:path4}), we obtain that the path (\ref{4to1:path3}) exists.
        
    This concludes our proof of the $\Sub(\psi)$-\accessibility of $\mathcal{F}^*$. 
    \end{proof}

\section{The finite frame property of $\NRAmn$, continued}\label{sec:ffp2}

This section is a continuation of the proof of our main theorem.
In particular, we prove the remaining case $n = 0$ of the implication (4 $\Rightarrow$ 1) of Theorem \ref{thm:ffp}. 
We note that $\NA{m}{0} = \NRA{m}{0}$ by Proposition \ref{prop:admissibility}.
So, we prove that for any $\psi \in \MF$, if $\NA{m}{0} \nvdash \psi$, then there exists a finite $\Sub(\psi)$-\accessible[1, 0] $\N$-frame on which $\psi$ is not valid. 
We first prove the case $m = 1$ since it is easier than the other cases.

\begin{lemma}[(4 $\Rightarrow$ 1) of Theorem \ref{thm:ffp} ($m = 1, n = 0$)] \label{lem:4to1-10}
  Take any $\psi \in \MF$ such that $\NA{1}{0} \nvdash \psi$,
  then there is a finite $\Sub(\psi)$-\accessible[1, 0] $\N$-model that falsifies $\psi$.
\end{lemma}
\begin{proof}
    Suppose that $\NA{1}{0} \nvdash \psi$, then $\set{ \negg \psi }$ is $\NA{1}{0}$-consistent.
    So there is $\psi$-maximal $\NA{1}{0}$-consistent set $w_\psi$ such that $\negg \psi \in w_\psi$.
    We let $M = (\set{ w_\psi }, \Rels, V)$ be an $\N$-model such that
    \begin{itemize}
        \item $w_\psi \Rel{\varphi} w_\psi$ iff $(\Box \varphi \notin w_\psi \tor \varphi \in w_\psi)$;
        \item $V(w_\psi, p) = 1$ iff $p \in w_\psi$.
    \end{itemize}
    The $\N$-model $M$ is clearly finite and $\Sub(\psi)$-\accessible[1, 0], so it suffices to show that $w_\psi \nVdash \psi$.
    \begin{claim*}
      For any $\rho \in \Sub(\psi)$, $\rho \in w_\psi$ if and only if $w_\psi \Vdash \rho$.
      \begin{subproof}[Proof of the claim]
        We use an induction on the construction of $\rho$.
        \begin{itemize}
          \item If $\rho = p \in \PropVar$, then $w_\psi \Vdash p \siff V(w_\psi, p) = 1 \siff p \in w_\psi$.
          \item If $\rho = \Box \rho'$, then we shall prove the claim in both directions.
                \begin{itemize}
                  \item[($\Rightarrow$)] Suppose that $w_\psi \nVdash \Box \rho'$, then $w_\psi \Rel{\rho'} w_\psi \nVdash \rho'$,
                    so $\rho' \notin w_\psi$ by the induction hypothesis.
                    This and $w_\psi \Rel{\rho'} w_\psi$ imply $\Box \rho' \notin w_\psi$ by the definition of $\Rel{\rho'}$.
                  \item[($\Leftarrow$)]
                    Suppose that $w_\psi \Vdash \Box \rho'$, then either $w_\psi \not \Rel{\rho'} w_\psi$ or $w_\psi \Rel{\rho'} w_\psi \Vdash \rho'$.
                    If $w_\psi \not \Rel{\rho'} w_\psi$, then $\Box \rho' \in w_\psi$ by the definition of $\Rel{\rho'}$.
                    Otherwise, $w_\psi \Vdash \rho'$ implies $\rho' \in w_\psi$ by the induction hypothesis.
                    This and $\NA{1}{0} \vdash \rho' \to \Box \rho'$ imply $\Box \rho' \in w_\psi$ by the maximal consistency of $w_\psi$.
                \end{itemize}
          \item The remaining cases ($\neg, \lor$) are trivial by the induction hypothesis. \qedhere
        \end{itemize}
      \end{subproof}
    \end{claim*}
    By the $\NA{1}{0}$-consistency of $w_\psi$, $\negg \psi \in w_\psi$ implies $\psi \notin w_\psi$.
    Therefore, $w_\psi \nVdash \psi$ by the claim.
  \end{proof}

However, the proof of (4 $\Rightarrow$ 1) of Theorem \ref{thm:ffp} for $m \ge 2$ and $n = 0$ is much more complicated.
Due to the nature of \accessibility[m,0], we construct a finite \accessible[m,0] $\N$-model falsifying $\psi$ depending on the shape of $\psi$,
which is described and analyzed with the notations introduced in the following definitions.

\begin{definition}
  Let $\varphi \in \MF$. 
  We define the natural numbers $\degs(\varphi)$ and $\degom(\varphi)$
  by induction on the construction of $\varphi$ as follows:

  \begin{itemize}
    \item $\degs(p) = \degom(p) : = 0$ for $p \in \PropVar$;
    \item $\degs(\neg \psi) : = \degs(\psi)$, \quad $\degom(\neg \psi) : = 0$;
    \item $\degs(\psi_1 \lor \psi_2) : = \max(\degs(\psi_1), \degs(\psi_2))$,
          \quad $\degom(\psi_1 \lor \psi_2) : = 0$;
    \item $\degs(\Box \psi) : = \begin{cases}
              \degs(\Box \rho) & \tif \psi = \Box \rho \ \text{for some}\ \rho,     \\
              \degs(\psi) + 1  & \mathrel{\text{otherwise}},
            \end{cases}$ \\
          $\degom(\Box \psi) : = \degom(\psi) + 1$.
  \end{itemize}

  Here, $\degs(\varphi)$ is a variant of modal degree in which consecutive boxes make up a single degree regardless of their number,
  and $\degom(\varphi)$ represents the number of the \emph{outer most} boxes. 
  So, $\degom(\varphi) > 0$ if and only if $\varphi$ is of the form $\Box \psi$. 
  Consider for example $\varphi = \Box (\Box p \lor \Box \Box q)$, then $\deg(\varphi) = 3$, $\degs(\varphi) = 2$, and $\degom(\varphi) = 1$.
\end{definition}

\begin{definition}
  Let $\varphi \in \MF$ and $0 \le d \le \degs(\varphi)$. We let
  \begin{align*}
    \Subdeg{\varphi, d}   & \defeq \set*{ \psi \in \Sub(\varphi); \degs(\psi) = d }, \\
    \SubdegLe{\varphi, d} & \defeq \bigcup_{0 \le i \le d} \Subdeg{\varphi, i}.
  \end{align*}
\end{definition}

  Consider for example $\varphi = \Box (\Box p \lor \Box \Box q)$,
  then $\Subdeg{\varphi, 0} = \set{p, q}$, $\Subdeg{\varphi, 1} = \set{(\Box p \lor \Box \Box q), \Box p, \Box \Box q, \Box q}$,
  and $\Subdeg{\varphi, 2} = \set{\Box (\Box p \lor \Box \Box q)}$.
  Note that $\SubdegLe{\varphi, \degs(\varphi)} = \Sub(\varphi)$.

\begin{definition}
  For any $\varphi \in \MF$ such that $\degom(\varphi) > 0$,
  it is easy to see that there is a unique $\psi \in \MF$ such that $\degs(\psi) = \degs(\varphi) - 1$, $\degom(\psi) = 0$, and $\varphi = \Box^{\degom(\varphi)} \psi$.
  We will call such a $\psi$ the \emph{root} of $\varphi$ and write it as $\Root(\varphi)$.
\end{definition}

In other words, $\Root(\varphi)$ is obtained from $\varphi$ with $\degom(\varphi) > 0$ by removing all the outer most boxes. 
We have $\degom(\Root(\varphi)) = 0$ and $\degs(\Root(\varphi)) = \degs(\varphi) - 1$. 
For example, $\Root(\Box(\Box p \lor \Box \Box q)) = \Box p \lor \Box \Box q$, and $\Root(\Box \Box q) = q$.

\begin{definition}
  Let $\varphi \in \MF$ and $1 \le d \le \degs(\varphi)$.
  We say $\lambda \in \Subdeg{\varphi, d}$ is a \emph{$d$-lord} on $\varphi$ if $\degom(\lambda) > 0$ and $\Box \lambda \notin \Subdeg{\varphi, d}$.
  We let $\Lord(\varphi, d)$ be the set of all $d$-lords on $\varphi$. 
  That is, 
  \[
    \Lord(\varphi, d) = \set*{ \lambda \in \Subdeg{\varphi, d}; \degom(\lambda) > 0, \Box \lambda \notin \Subdeg{\varphi, d} } .
  \]
\end{definition}

Consider for example $\varphi = \Box (\Box p \lor \Box \Box q)$, then $\Lord(\varphi, 1) = \set{ \Box p, \Box \Box q }$
  and $\Lord(\varphi, 2) = \set{ \varphi }$.

\begin{definition}
    Let $\lambda \in \MF$ be such that $\degom(\lambda) > 0$. 
    We say $\rho$ is a \emph{liege} of $\lambda$ if $\degom(\rho) > 0$ and $\lambda = \Box^i \rho$ for some $i > 0$. 
  We let $\Liege(\lambda)$ be the set of all lieges of $\lambda$. 
\end{definition}
   
  If $\rho$ is a liege of $\lambda$, then $\degs(\rho) = \degs(\lambda)$, $\Root(\rho) = \Root(\lambda)$, and $0 < \degom(\rho) < \degom(\lambda)$. 
  We have 
 \[
    \Liege(\lambda) = \set*{ \Box^i \Root(\lambda); 0 < i < \degom(\lambda) }.
\]
  For example, $\Liege(\Box p) = \emptyset$, $\Liege(\Box \Box q) = \set{ \Box q }$, and $\Liege(\Box^3 r) = \set{ \Box \Box r, \Box r }$.

  We shall note a couple of properties on $d$-lords and their lieges:

  \begin{itemize}
    \item A root is neither a lord nor a liege.
    \item Let $\lambda$ be a $d$-lord, then $\Subdeg{\lambda, d} = \set{\lambda} \cup \Liege(\lambda)$.
    \item For any $d$-lords $\lambda_1$ and $\lambda_2$ on $\varphi$, $\Root(\lambda_1) = \Root(\lambda_2)$ if and only if $\lambda_1 = \lambda_2$.
          Therefore, $\lambda_1 \ne \lambda_2$ implies $\Liege(\lambda_1) \cap \Liege(\lambda_2) = \emptyset$.
  \end{itemize}

Now we shall construct a finite \accessible[m,0] $\N$-model falsifying $\psi$ by induction on $\degs(\psi)$;
we start with a model with a single world and no accessibility relations to falsify propositional formulae,
then we repeatedly add only the necessary worlds and relations to falsify formulae with more modal degree.
We introduce a notation for finite $\N$-models so that we can easily state that a model contains only the desired relations.

\newcommand*{\angled}[1]{\langle #1 \rangle}

\begin{definition}
  For any finite set $W$ and any finite $R \subseteq W \times \MF \times W$,
  we write $\angled{W, R}$ to mean the finite $\N$-frame $(W, \Rels)$
  where for any $\varphi \in \MF$ and $w, w' \in W$, 
  \[
    w \Rel{\varphi} w' \iff (w, \varphi, w') \in R.
  \]
  For the sake of brevity, we write $(w \Rel{\varphi} w') \in R$ to mean $(w, \varphi, w') \in R$.
  We also write $\angled{W, R, V}$ to mean a finite $\N$-model in the same manner.
\end{definition}

  For example, $\angled{W, \emptyset}$ is the finite $\N$-frame $(W, \Rels)$ where $w \not \Rel{\varphi}^* w'$ for all $\varphi \in \MF$ and $w, w' \in W$.

We are ready to prove our theorem on the existence of falsifying $\N$-models. 

\begin{theorem} \label{thm:m0-restricted-canonical-model}
  Let $m \ge 2$. Take any $\varphi \in \MF$,
  any $\psi \in \Sub(\varphi)$ such that $\NA{m}{0} \nvdash \psi$,
  and any $\varphi$-maximal $\NA{m}{0}$-consistent set $w_\psi$ such that $\negg \psi \in w_\psi$.
  For each $0 \le d \le \degs(\psi)$, there is a finite $\N$-model $M_{\psi, d} = \angled{W_{\psi, d}, R_{\psi, d}, V_{\psi, d}}$
  such that $w_\psi \in W_{\psi, d}$ and:
  \begin{enumerate}[label=(\alph*)]
    \item \label{m0-rcm:a} $w_\psi \Vdash \rho$ if and only if $\rho \in w_\psi$ for any $\rho \in \SubdegLe{\psi, d}$;
    \item \label{m0-rcm:b} $\angled{W_{\psi, d}, R_{\psi, d}}$ is $\Sub(\varphi)$-\accessible[m,0];
    \item \label{m0-rcm:c} if $d \ge 1$, then for any $(w \Rel{\rho} w') \in R_{\psi, d}$, either $\degs(\rho) < d$
          or $\rho$ is a liege of some $d$-lord on $\psi$, that is, $\rho \in \bigcup_{\lambda \in \Lord(\psi, d)} \Liege(\lambda)$.
  \end{enumerate}
\end{theorem}
  \begin{proof}
    We prove the theorem by induction on $\degs(\psi)$ (we will refer to this induction as \hypertarget{m0-rcm:ind-1}{\emph{Induction 1}}).

    Suppose first that $\degs(\psi) = 0$ and take any $0 \le d \le \degs(\psi)$, then $d = 0$.
    We let $M_{\psi, 0} = \angled{\set{w_\psi}, \emptyset, V}$,
    where $V(w_\psi, p) = 1 \defsiff p \in w_\psi$ for every $p \in \PropVar$.
    This model is clearly finite, and \ref{m0-rcm:b} and \ref{m0-rcm:c} trivially hold since $R_{\psi, 0} = \emptyset$ and $d = 0$.
    For \ref{m0-rcm:a}, take any $\rho \in \SubdegLe{\psi, 0}$,
    then $d = 0$ implies that $\degs(\rho) = 0$. 
    It easily follows from the $\varphi$-maximal $\NA{m}{0}$-consistency of $w_\psi$ that $w_\psi \Vdash \rho$ if and only if $\rho \in w_\psi$.

    Suppose next that $\degs(\psi) > 0$ and the theorem holds for any $\psi' \in \Sub(\varphi)$ such that $0 \le \degs(\psi') < \degs(\psi)$. 
    We prove the theorem for $\psi$ by induction on $d$ (\hypertarget{m0-rcm:ind-2}{\emph{Induction 2}}).

    Suppose that $d = 0$, then we let $M_{\psi, 0} = \angled{\set{w_\psi}, \emptyset, V}$,
    where $V(w_\psi, p) = 1 \defsiff p \in w_\psi$ for every $p \in \PropVar$.
    It is verified that $M_{\psi, 0}$ satisfies \ref{m0-rcm:a}, \ref{m0-rcm:b} and \ref{m0-rcm:c} in the same way as above.

    Now suppose that $d > 0$ and the theorem holds for $d-1$.

    Let $\set{ \lambda_1, \lambda_2, \ldots \lambda_{\left|\Lord(\psi, d)\right|} } = \Lord(\psi, d)$.
    We shall inductively define finite models
    $M_{\psi, d, 0},\ M_{\psi, d, 1},\ \ldots,\ M_{\psi, d, \left|\Lord(\psi, d)\right|}$
    such that, for each $i \leq |\Lord(\psi, d)|$, we have $M_{\psi, d, i} = \angled{W_{\psi, d, i}, R_{\psi, d, i}, V_{\psi, d, i}}$, $w_\psi \in W_{\psi, d, i}$ and:

    \begin{enumerate}[label=(${\text{\alph*}}^\ast$)]
      \item \label{m0-rcm:astar} $w_\psi \Vdash \rho$ if and only if $\rho \in w_\psi$ for any $\rho$ such that
            either $\rho \in \SubdegLe{\psi, d-1}$ or $\rho \in \bigcup_{1 \le x \le i} \Subdeg{\lambda_x, d}$;
      \item \label{m0-rcm:bstar} $\angled{W_{\psi, d, i}, R_{\psi, d, i}}$ is $\Sub(\varphi)$-\accessible[m,0];
      \item \label{m0-rcm:cstar} For any $(w \Rel{\rho} w') \in R_{\psi, d, i}$, $\rho$ satisfies one of the following:
            \begin{itemize}
              \item $\degs(\rho) < d-1$;
              \item $\degs(\rho) = d-1$ and $\rho$ is a liege of some $(d-1)$-lord on $\psi$;
              \item $\degs(\rho) = d-1$ and $\rho$ is a root of one of $\lambda_1, \ldots, \lambda_i$;
              \item $\degs(\rho) = d$ and $\rho$ is a liege of one of $\lambda_1, \ldots, \lambda_i$.
            \end{itemize}
    \end{enumerate}

    Then let $M_{\psi, d} : = M_{\psi, d, \left|\Lord(\psi, d)\right|}$. Assuming that $M_{\psi, d}$ is defined, we first check that
    $M_{\psi, d}$ satisfies \ref{m0-rcm:a}, \ref{m0-rcm:b} and \ref{m0-rcm:c}. 
    It is clear that \ref{m0-rcm:bstar} and \ref{m0-rcm:cstar} imply \ref{m0-rcm:b} and \ref{m0-rcm:c} respectively
    by letting $i = \left|\Lord(\psi, d)\right|$. 
    For \ref{m0-rcm:a}, take any $\rho \in \SubdegLe{\psi, d}$. If $\degs(\rho) < d$,
    then $\rho \in \SubdegLe{\psi, d-1}$, so $w_\psi \Vdash \rho$ if and only if $\rho \in w_\psi$ by \ref{m0-rcm:astar}.
    Otherwise, $\degs(\rho) = d$. 
    We prove the equivalence \ref{m0-rcm:a} for $\rho$ by induction on the construction of $\rho$.
    \begin{itemize}
      \item Since $d > 0$, it is impossible that $\rho \in \PropVar$.
      \item If $\rho$ is of the form $\Box \rho'$, then $\degom(\rho) > 0$, so we have 
      \[
        \rho \in \bigcup_{1 \le x \le \left|\Lord(\psi, d)\right|} \Subdeg{\lambda_x, d}.
      \]
            Therefore, $w_\psi \Vdash \rho$ if and only if $\rho \in w_\psi$ by \ref{m0-rcm:astar}, letting $i = \left|\Lord(\psi, d)\right|$.
      \item The remaining cases ($\neg, \lor$) are trivial by the induction hypothesis.
    \end{itemize}

    Now we construct $M_{\psi, d, i}$ for $0 \le i \le \left|\Lord(\psi, d)\right|$ by induction on $i$.

    We first construct $M_{\psi, d, 0}$.
    By the hypothesis of \hyperlink{m0-rcm:ind-2}{\emph{Induction 2}}, the model $M_{\psi, d-1}$ exists and satisfies \ref{m0-rcm:a}, \ref{m0-rcm:b} and \ref{m0-rcm:c}.
    We let $M_{\psi, d, 0} : = M_{\psi, d-1}$. 
    Then $M_{\psi, d, 0}$ satisfies \ref{m0-rcm:astar}, \ref{m0-rcm:bstar} and \ref{m0-rcm:cstar} because
    $\bigcup_{1 \le x \le 0} \Subdeg{\lambda_x, d} = \emptyset$,
    and for any $(w \Rel{\rho} w') \in R_{\psi, d-1}$, either $\degs(\rho) < d-1$ or
    $\rho$ is a liege of some $(d-1)$-lord on $\psi$.

    Now suppose that $i > 0$ and $M_{\psi, d, i-1}$ is defined.
    Let $\chi : = \Root(\lambda_i)$, then $\degs(\chi) = d-1$, $\degom(\chi) = 0$, and $\lambda_i = \Box^{\degom(\lambda_i)} \chi$.
    For $0 \le j < m$, we let
    \begin{equation*}
      s_j : = \max \left( \set{0} \cup \set*{ 1 \le k \le \degom(\lambda_i); k \equiv j \Mod m \;\&\; \Box^k \chi \notin w_\psi } \right).
    \end{equation*}

    We first make a claim on $s_j$:

    \begin{claim*}
      For any $1 \le k \le \degom(\lambda_i)$ such that $k \equiv j \Mod m$,
      $$\Box^k \chi \in w_\psi \iff k > s_j.$$
      \begin{subproof}[Proof of the claim]
        ($\Leftarrow$)
        Trivial by the definition of $s_j$.
        ($\Rightarrow$)
        Suppose towards a contradiction that $\Box^k \chi \in w_\psi$ but $k \le s_j$.
        Since $1 \leq k$, we have $s_j \neq 0$, and so $s_j \equiv j \equiv k \Mod m$, which implies $\NA{m}{0} \vdash \Box^k \chi \to \Box^{s_j} \chi$. 
        We have $\Box^{s_j} \chi \notin w_\psi$ by the definition of $s_j$, 
        so $\Box^k \chi \notin w_\psi$ by the maximal consistency of $w_\psi$,
        which is a contradiction.
      \end{subproof}
    \end{claim*}

    Now we shall construct $M_{\psi, d, i}$. We distinguish the following cases:

    \subparagraph*{Case 1:} Suppose that $s_0 = s_1 = \cdots = s_{m-1} = 0$, then this implies $\Subdeg{\lambda_i, d} = \set{ \Box \chi, \Box^2 \chi, \ldots, \lambda_i } \subseteq w_\psi$.
    We let $M_{\psi, d, i} : = M_{\psi, d, i-1}$. 
    Then, \ref{m0-rcm:bstar} and \ref{m0-rcm:cstar} trivially hold since $R_{\psi, d, i} = R_{\psi, d, i-1}$.
    For \ref{m0-rcm:astar}, by \ref{m0-rcm:astar} of $M_{\psi, d, i-1}$, it suffices to show that $w_\psi \Vdash \rho$ if and only if $\rho \in w_\psi$ for every $\rho \in \Subdeg{\lambda_i, d}$.
    Here, $\Subdeg{\lambda_i, d} \subseteq w_\psi$ implies that we only need to show that $w_\psi \Vdash \rho$ for every $\rho \in \Subdeg{\lambda_i, d}$.
    Let $\rho = \Box \rho'$, then $\rho'$ is a liege or the root of $\lambda_i$, so $\rho'$ cannot be a liege or the root of $\lambda_x$ for any $0 \le x \le i-1$.
    Here, $(w \Rel{\rho'} w') \notin R_{\psi, d, i} = R_{\psi, d, i-1}$ by \ref{m0-rcm:cstar} of $M_{\psi, d, i-1}$.
    Therefore, $w_\psi \Vdash \Box \rho'$.

    \subparagraph*{Case 2:} Now suppose that there is at least one $j$ such that $s_j \ne 0$,
    then $\Box^{s_j} \chi \notin w_\psi$ implies $\NA{m}{0} \nvdash \Box^{s_j} \chi$ by the maximal consistency of $w_\psi$,
    which implies $\NA{m}{0} \nvdash \chi$ by applying \textsc{Nec} for $s_j$ times.
    This implies that $\set{ \negg \chi }$ is $\NA{m}{0}$-consistent, so there is a $\varphi$-maximal $\NA{m}{0}$-consistent
    set $w_\chi$ such that $\negg \chi \in w_\chi$.
    Here, $\degs(\chi) = d-1 < d = \degs(\psi)$, so by the hypothesis of \hyperlink{m0-rcm:ind-1}{\emph{Induction 1}},
    there is a finite $\N$-model $M_{\chi, d-1} = \angled{W_{\chi, d-1}, R_{\chi, d-1}, V_{\chi, d-1}}$ satisfying \ref{m0-rcm:a}, \ref{m0-rcm:b} and \ref{m0-rcm:c} and $w_\chi \in W_{\chi, d-1}$.

    For each $1 \le j < m$, we let $M_{\chi, d-1, j} = \angled{W_{\chi, d-1, j}, R_{\chi, d-1, j}, V_{\chi, d-1, j}}$ be a copy of $M_{\chi, d-1}$ which is disjoint from $M_{\psi, d, i-1}$ and from $M_{\chi, d-1, x}$ for every $1 \le x < j$.
    We also let $w_{\chi, j} \in W_{\chi, d-1, j}$ be the copy of $w_\chi$, then $w_{\chi, j} \Vdash_{M_{\chi, d-1, j}} \negg \chi$.

    We define a sequence $\{w_l^*\}$ as follows: 
    \begin{equation*}
      w_l^* := \begin{cases}
        w_\psi      & \tif l \equiv 0 \Mod m, \\
        w_{\chi, j} & \tif l \equiv j \Mod m.
      \end{cases}
    \end{equation*}
    It is easy to see that $w_{l + m}^* = w_l^*$ for any $l \ge 0$.

    We define the model $M_{\psi, d, i} = \angled{W_{\psi, d, i}, R_{\psi, d, i}, V_{\psi, d, i}}$ as follows:

    \begin{itemize}
      \item $W_{\psi, d, i} : = W_{\psi, d, i-1} \:\cup\: \bigcup_{1 \le j < m} W_{\chi, d-1, j}$;
      \item $R_{\psi, d, i} : = R_{\psi, d, i-1} \:\cup\: \bigcup_{1 \le j < m} R_{\chi, d-1, j} \:\cup\: \bigcup_{0 \le j < m} R_j$, where \\
            \begin{equation*}
              R_j = \set*{ w_l^* \Rel{\Box^{s_j - l - 1} \chi} w_{l+1}^*; 0 \le l < s_j };
            \end{equation*}
      \item $V_{\psi, d, i}(w, p) : = \begin{cases}
                V_{\psi, d, i-1}(w, p) & \tif w \in W_{\psi, d, i-1}, \\
                V_{\chi, d-1, j}(w, p) & \tif w \in W_{\chi, d-1, j}.
              \end{cases}$
    \end{itemize}

\begin{figure}[ht]
\centering
\begin{tikzpicture}
\draw (0, 0) circle [radius=1];
\draw (0, 1) node[above]{$M_{\psi, d, i-1}$};
\fill (0, 0) circle (1pt);
\draw (0, 0) node[above]{$w_\psi$};

\draw (2.5, 0) circle [radius=1];
\draw (2.5, 1) node[above]{$M_{\chi, d-1, 1}$};
\fill (2.5, 0) circle (1pt);
\draw (2.5, 0) node[above]{$w_{\chi, 1}$};

\draw (5, 0) node {$\cdots$};

\draw (7.5, 0) circle [radius=1];
\draw (7.5, 1) node[above]{$M_{\chi, d-1, m-1}$};
\fill (7.5, 0) circle (1pt);
\draw (7.5, 0) node[above]{$w_{\chi, m-1}$};

\draw [thick, ->, dashed] (0, 0) to [out=-45,in=225](2.5,0);
\draw [thick, ->, dashed] (2.5, 0) to [out=-45,in=180](3.75,-0.5);
\draw [thick, ->, dashed] (6.25, -0.5) to [out=0,in=225](7.5,0);
\draw [thick, ->, dashed] (7.5, 0) to [out=235,in=305](0,0);
\end{tikzpicture}
\caption{The model $M_{\psi, d, i}$}
    \label{fig3}
\end{figure}

    The model $M_{\psi, d, i}$ is visualized in Figure \ref{fig3}. 
    Here, $M_{\psi, d, i}$ is finite and satisfies \ref{m0-rcm:cstar} by the definition of $R_{\psi, d, i}$. 
    For each $0 \le j < m$ such that $s_j > 0$, by the definition of $R_j$, it is clear that the following path $P_j$ exists:
    \begin{equation*}
      P_j \ \defeq\  w_0^* \Rel{\Box^{s_j - 1} \chi} w_1^* \Rel{\Box^{s_j - 2} \chi} \ \cdots\ \Rel{\Box \chi} w_{s_j -1}^* \Rel{\chi} w_{s_j}^*.
    \end{equation*}

    Now we prove that $M_{\psi, d, i}$ satisfies \ref{m0-rcm:astar} and \ref{m0-rcm:bstar}.

    \paragraph*{\ref{m0-rcm:astar} of $M_{\psi, d, i}$:}
    Take any $\rho \in \SubdegLe{\psi, d-1} \cup \bigcup_{1 \le x \le i} \Subdeg{\lambda_x, d}$,
    then we shall show that $w_\psi \Vdash \rho$ if and only if $\rho \in w_\psi$.
    We prove this equivalence by induction on the construction of $\rho$.
    We will only show the case of $\rho = \Box \sigma$, as the other cases are trivial by using the induction hypothesis.

    If $\sigma$ is not a liege of $\lambda_i$,
    then by the definition of $R_{\psi, d, i}$,
    for every $w \in W_{\psi, d, i}$, $(w_\psi \Rel{\sigma} w) \in R_{\psi, d, i}$ iff $(w_\psi \Rel{\sigma} w) \in R_{\psi, d, i-1}$.
    This and the definition of $V_{\psi, d, i}$ imply
    $w_\psi \Vdash \Box \sigma$ iff $w_\psi \Vdash_{M_{\psi, d, i-1}} \Box \sigma$ iff $\Box \sigma \in w_\psi$.

    Now we suppose that $\sigma$ is a liege of $\lambda_i$.
    Let $\Box^k \chi = \Box \sigma$, then $1 \le k \le \degom(\lambda_i)$, and
    there is $0 \le j < m$ such that $k \equiv j \Mod m$.
    Here, $\Box^k \chi \in w_\psi$ iff $k > s_j$ by the claim.
    So we shall show that the equivalence 
    \[
        w_\psi = w_0^* \Vdash \Box^k \chi \iff k > s_j
    \]
    holds in both directions.

    \subparagraph*{($\Rightarrow$)}
    If $(0 < ) k \le s_j$, then the path $P_j$ exists,
    and $k \equiv j \equiv s_j \Mod m$ implies that:
    \begin{equation*}
      w_0^* = w_{s_j - k}^* \Rel{\Box^{k-1} \chi} w_{s_j - k + 1}^* \Rel{\Box^{k-2} \chi} \ \cdots\ \Rel{\Box \chi} w_{s_j - 1}^* \Rel{\chi} w_{s_j}^* = w_j^*.
    \end{equation*}
    If $j \ne 0$, then $w_j^* = w_{\chi, j}$.
    Since $w_{\chi, j} \Vdash_{M_{\chi, d-1, j}} \negg \chi$,
    it is easy to show that $w_{\chi, j} \Vdash \negg \chi$,
    so $w_0^* \nVdash \Box^k \chi$.
    Otherwise, $j = 0$ and $w_j^* = w_0^* = w_\psi$.
    Here, $k \equiv j \equiv 0 \Mod m$ implies $\NA{m}{0} \vdash \chi \to \Box^k \chi$,
    and $k \le s_0$ implies $\Box^k \chi \notin w_\psi$,
    so $\chi \notin w_\psi$ by the maximal consistency of $w_\psi$.
    Therefore, $w_\psi \nVdash \chi$ by the induction hypothesis,
    which implies $w_\psi \nVdash \Box^k \chi$.

    \subparagraph*{($\Leftarrow$)}
    If $k > s_j$, then $k - 1 \ge s_j$, so by the definition of $R_{\psi, d, i}$,
    there is no $w \in W_{\psi, d, i}$ such that $w_\psi \Rel{\Box^{k-1} \chi} w$,
    which implies $w_\psi \Vdash \Box^k \chi$.

    \paragraph*{\ref{m0-rcm:bstar} of $M_{\psi, d, i}$:}

    Take any path of length $m$:
    $$
      w_0 \Rel{\Box^{m-1} \rho} w_1 \Rel{\Box^{m-2} \rho} w_2 \ \cdots\ w_{m-2} \Rel{\Box \rho} w_{m-1} \Rel{\rho} w_m,
    $$
    and we shall show that $w_0 = w_m$.
    Let $R = \set{ (w_0 \Rel{\Box^{m-1} \rho} w_1), \ldots, (w_{m-1} \Rel{\rho} w_m) }$.

    If either $\rho = \chi$ or $\rho$ is a liege of $\lambda_i$,
    then $R \subseteq \bigcup_{0 \le j < m} R_j$ by \ref{m0-rcm:cstar} of $M_{\psi, d, i-1}$ and \ref{m0-rcm:c} of $M_{\chi, d-1}$.
    Then it is easy to show that this path is a sub-path of some $P_j$, and there is $0 \le l < s_j$ such that $w_0 = w_l^*, w_1 = w_{l+1}^*, \ldots, w_m = w_{l + m}^*$ and:
    \begin{equation*}
      w_l^* \Rel{\Box^{s_j - l - 1} \chi} w_{l+1}^* \Rel{\Box^{s_j - l - 2} \chi} \cdots \Rel{\Box^{s_j - l - m} \chi} w_{l + m}^*.
    \end{equation*}
    Here, $w_0 = w_l^* = w_{l+m}^* = w_m$.

    Otherwise, $R \cap \bigcup_{0 \le j < m} R_j = \emptyset$,
    then by the definition of $M_{\psi, d, i}$ and the disjointness of the models,
    either $R \subseteq R_{\psi, d, i-1}$
    or $R \subseteq R_{\chi, d-1, k}$ for some $1 \le j < m$.
    Therefore, $w_0 = w_m$ by \ref{m0-rcm:bstar} of $M_{\psi, d, i-1}$ and \ref{m0-rcm:b} of $M_{\chi, d-1}$ respectively.
  \end{proof}

\begin{lemma}[(4 $\Rightarrow$ 1) of Theorem \ref{thm:ffp} ($m \ge 2, n = 0$)]
  Let $m \ge 2$. Take any $\varphi \in \MF$ such that $\NA{m}{0} \nvdash \varphi$,
  then there is a finite $\Sub(\varphi)$-\accessible[m, 0] $\N$-model that falsifies $\varphi$.
\end{lemma}
  \begin{proof}
    Suppose that $\NA{m}{0} \nvdash \varphi$, then $\set{ \negg \varphi }$ is $\NA{m}{0}$-consistent.
    So there is a $\varphi$-maximal $\NA{m}{0}$-consistent set $w_\varphi$ such that $\negg \varphi \in w_\varphi$.
    Then by Theorem \ref{thm:m0-restricted-canonical-model},
    there is a finite $\Sub(\varphi)$-\accessible[m, 0] model $M_{\varphi, \degs(\varphi)} = \angled{W, R, V}$ such that $w_\varphi \in W$ and
    $w_\varphi \Vdash \rho$ if and only if $\rho \in w_\varphi$ for every $\rho \in \SubdegLe{\varphi, \degs(\varphi)} = \Sub(\varphi)$.
    Here, $\varphi \notin w_\varphi$ by the $\varphi$-maximal $\NA{m}{0}$-consistency, so $w_\varphi \nVdash \varphi$.
  \end{proof}

This concludes the proof of (4 $\Rightarrow$ 1) of Theorem \ref{thm:ffp}.
We obtain the following corollary.

\begin{corollary}
  For each $m,n \in \mathbb{N}$, the logic $\NRAmn$ is decidable.
\end{corollary}

\section{Future Work}\label{sec:futurework}

In this paper, we proved the finite frame property of the the extensions of $\N$ with the axiom scheme $\Amn: \Box^n \varphi \to \Box^m \varphi$ and the rule $\RosBox$.
The obvious next work would be to prove the same with multiple axioms schemata of the form $\Amn$.
For any $X \subseteq \mathbb{N}^2$, let $\mathbf{NA}_X$ be the logic obtained from $\N$ by adding the axiom schemata $\Amn$ for $(m,n) \in X$. 
We say that an $\N$-frame $\mathcal{F}$ is $X$-accessible if it is \accessible for all $(m, n) \in X$. 
It follows from Theorem \ref{thm:soundness}, the logic $\mathbf{NA_X}$ is sound with respect to $X$-accessible $\N$-frames. 
Also let $\mathbf{NA}_X^+$ and $\mathbf{NRA}_X$ be the logics obtained from $\mathbf{NA}_X$ by adding the rules $\RosBox$ and \textsc{Ros}, respectively.  
Then, it is shown that the logic $\mathbf{NRA}_X$ is sound with respect to serial $X$-accessible $\N$-frames. 

\begin{problem}
Let $X \subseteq \mathbb{N}^2$. 
\begin{enumerate}
    \item Is $\mathbf{NA}_X^+$ sound with respect to the class of all $X$-accessible $\N$-frames?
    \item Is either $\mathbf{NA}_X$ or $\mathbf{NA}_X^+$ complete with respect to the class of all $X$-accessible $\N$-frames?
    \item Is $\mathbf{NRA}_X$ complete with respect to the class of all serial $X$-accessible $\N$-frames?
\end{enumerate}
\end{problem}

In particular, the proof of Lemma \ref{lem:3to4} ((3 $\Rightarrow$ 4) of Theorem \ref{thm:ffp}) does not seem to be easily generalized; Some $(m_1, n_1) \in X$ may fall into the case \ref{3to4:relstar2} while some other $(m_2, n_2) \in X$ may fall into \ref{3to4:relstar3} in the proof of Lemma \ref{lem:3to4}. 
Some modification on the relation $\Rel{\varphi}^*$ would be needed to obtain a general result.


In the context of provability logic, one may wonder the completeness of $\N + \mathrm{GL}: \Box \left(\Box \varphi \to \varphi \right) \to \Box \varphi$,
which is not of the form $\Amn$ and thus not covered by this paper.
Many of the $\N$ counterpart of the famous extensions of $\K$ are still left uninvestigated. 

We have only investigated the completeness and the finite frame property of $\NRAmn$,
so it may also be interesting to investigate other good properties such as compactness and the interpolation property.

\section*{Acknowledgements}

The authors would like to thank Haruka Kogure for pointing out some errors in our proof presented in the initial version of this paper. The authors also thank the anonymous referees for giving helpful advice on improving this paper. The first author was supported by JSPS KAKENHI Grant Numbers JP19K14586 and JP23K03200.

\bibliographystyle{plain}
\bibliography{refs}

\end{document}